\documentclass[11pt]{amsart}
\usepackage{amsmath,amscd,amssymb, latexsym, amsthm}
\newcommand{\cw}{\curlywedge}
\newcommand{\A}{{\mathcal{A}}}

\newcommand{\I}{{\mathcal{I}}}
\newcommand{\chivar}{$\chi$-invariant}
\newcommand{\del}{\partial}
\newcommand{\C}{{\mathbb C}}

\newtheorem{cor}{Corollary}
\newtheorem{prop}{Proposition}
\newtheorem{lemma}{Lemma}

\title{Semiinvariants of Finite Reflection Groups}
\author{Anne V. Shepler}

\address{Department of Mathematics,
       University of California at San Diego,
       La Jolla, California, 92093-0112}
\email{ ashepler@euclid.ucsd.edu}

\subjclass{51F15 (Reflection groups, reflection geometries), 
          52B30 (Arrangements of hyperplanes), 20 (Group theory).}

\begin{document}

\begin{abstract}
Let $G$ be a finite group of complex $n \times  n$ unitary matrices generated
by reflections acting on $\C^{\hspace{.1ex}n}$.  Let $R$
be the ring of invariant polynomials, and $\chi$ be a multiplicative
character of $G$.  Let $\Omega^\chi$ be the $R$-module 
of \chivar{} differential forms.  We define a multiplication 
in $\Omega^\chi$ and show that under this multiplication 
$\Omega^\chi$ has 
an {\em exterior algebra} structure.  We also show how to 
extend the results to vector
fields, and exhibit a relationship between
\chivar{} forms and logarithmic forms. 
\end{abstract}
\maketitle
\section{Introduction}
In 1989, P. Doyle and C. McMullen \cite{McM-Doy} solved the fifth degree polynomial with
a highly symmetrical dynamical system which preserved the Galois group $A_5$.
In 1997, S. Crass and P. Doyle \cite{Cra-Doy} solved the sixth degree polynomial
by again finding a dynamical system with special symmetry---this time $A_6$
symmetry.  Each dynamical system was formed by iterating a map that
was equivariant under the projective action of the group.  Such maps
correspond naturally to semiinvariant differential forms.
Because almost nothing was known about these forms, constructing the necessary
dynamical systems was a difficult step in both cases. 

We introduce here a general theory of semiinvariants.  
Specifically, we show that for any finite unitary reflection
group $G$ and multiplicative character $\chi$ of $G$, the module of $\chi$-invariant differential 
forms has a natural multiplication which
turns the module into an {\em exterior algebra}.  This exterior algebra
structure allows us to understand completely the forms that give rise to
highly symmetrical dynamical systems, and gives us tools to compute these
forms explicitly.  We also show how to extend these results to vector fields (or
{\em derivations}), and observe the relationship between semiinvariants
and logarithmic forms.  

The theory presented here builds on work by R. Stanley, who
characterized the module of $\chi$-invariant polynomials in 1977 \cite{Stan}.  It also
builds on more recent work by Orlik, Saito, Solomon, Terao
and others on invariant derivations and the theory of hyperplane
arrangements (see \cite{OrTer}, 
Chapter 6).  Note that
$\det$-invariant forms have received attention under the name of 
{\em anti-invariant forms} 
in the context of Coxeter groups (see e.g. \cite{double}).

\section{notation}
Let $G$ be a finite group of complex $n \times  n$ unitary matrices generated
by reflections acting on $V := \C^{\hspace{.1ex}n}$.  Recall that
a unitary matrix is a reflection if it has finite order 
and fixes a hyperplane pointwise in $V$.  Let $S := \C[x_1, \ldots, x_n]$ 
be the ring
of polynomials of $V$.    
Let $f_1, \ldots, f_n \in S$ be basic invariants, 
and $R = \C[f_1, \ldots, f_n]$ be the ring of invariant polynomials.
Let $\chi$ be a multiplicative
character of $G$. 
Denote the module of differential $p$-forms on $V$ by 
\begin{eqnarray*}
    \Omega^p & := & \bigoplus_{1 \leq i_1 < \ldots < i_p \leq n}
    S dx_{i_1} \wedge \ldots \wedge dx_{i_p}  \\
    &\simeq & S \otimes {\textstyle\bigwedge}^p V^* .
\end{eqnarray*}
The group $G$ acts contragradiently on $V^*$ and $S$, and
$\Omega^p$ is a $\C[G]$-module.
Define the $R$-module of {\em \chivar{}} differential $p$-forms as
$$
  \begin{array}{rcl}
     (\Omega^p)^{\chi} & := & \left\{\omega \in \Omega^p: g \omega = 
          \chi(g) \omega \text{ for all } g \in G \right\}.
  \end{array}
$$
Let
  \begin{eqnarray*}
     \Omega^{\chi} & := & \bigoplus_{0 \leq p} (\Omega^p)^{\chi}.
  \end{eqnarray*}

It is convenient to define 
$ \I^p $ as the set of multiindices of $\{ 1,...,n \}$ of length $p$:
$$
  \I^p := \left\{ I = \{ I_1, ..., I_p \} : 1 \leq I_1 < \ldots < I_p \leq n 
   \right\}.
$$
For a multiindex $I$, let $I^c$ denote the complementary index.
Denote the volume
form on $V$ by $vol :=  dx_1 \wedge \ldots \wedge dx_n$.  
If $f$ and $g$ are differential forms, we write $f \doteq g$
if $f = cg$ for some $c$ in $\C^*$.

We recall some facts and notation from \underline{Arrangements of
Hyperplanes} (\cite{OrTer}, p. 228).  Let
$\A$ be the hyperplane arrangement defined by $G$.  For each
$H \in \A$, define $\alpha_H \in S$ by 
ker$(\alpha_H) = H$.   
Fix some $H \in \A$,
and let $G_H$ be the cyclic subgroup of elements in $G$ that
fix $H$ pointwise.  Let $s_H$ be a generator of $G_H$ and 
let $o(s_H)$ be the order of $s_H$.  Define
$a_H(\chi)$ as the least integer satisfying
$0 \leq a_H(\chi) < o(s_H)$ and
$\chi(s_H)= \det(s_H)^{-a_H(\chi)}$.
Let
$$
   Q_\chi = \prod_{H \in \A} \alpha_H^{{a_H}(\chi)}.
$$
The polynomial $Q_\chi$ is uniquely determined, upto a nonzero scaler 
multiple, by the group $G$.  

R. Stanley ~\cite{Stan} proved that   
$   (\Omega^0)^{\chi} = 
               R \thinspace Q_{\chi}$,
and since $vol$ is $(\det^{-1})$-invariant, 
it follows that  
\begin{equation}
     \label{invn-forms}
 \tag{$*$}    (\Omega^n)^\chi = 
               R \thinspace Q_{\chi\cdot\det}\thinspace vol.
\end{equation}
R. Steinberg ~\cite{Steinberg} proved  that 
$
    Q_{\det} = \prod_{H \in \A} \alpha_H^{o(s_H)-1}
$
is the determinant of the Jacobian matrix 
$\left\{ \frac{\displaystyle \del } 
{\displaystyle \del x_i} f_j \right\} $, upto a nonzero scalar multiple.  
Note also 
that
$
  Q_{\det^{-1}} = \prod_{H \in \A} \alpha_H
$
(\cite{OrTer}, p. 229).

\section{$\chi$-wedging}

The next lemma will be used to show that $Q_\chi$ divides the exterior
product of any two \chivar{} forms.

\begin{lemma}
\label{lemma:degx1}
  Suppose that $\mu$ 
  is a \chivar{} $p$-form.  Fix a hyperplane $H \in \A$, and 
  let $a= a_H(\chi)$. 
  Choose coordinates in which $x_1 = \alpha_H$ and $s_H$ is diagonal.
  If   
  $$
    \mu =  \displaystyle \sum_{I \in \I^p} \mu_I \thinspace dx_{I_1} 
         \wedge \hdots \wedge dx_{I_p}
  $$ 
  in these coordinates, then  
  $x_1^{a-1}$ divides $\mu_I$ whenever $I_1 = 1$ and 
  $ x_1^{a}$ divides $\mu_I$  whenever $I_1 \neq 1$, 
   for each $I = \{I_1, \ldots, I_p\} \in \I^p$. 
 \end{lemma}


\begin{proof}
Let $s = s_H$ and $\rho$ be the 
determinant of $s$. Then 
$$
s = \left(
\begin{matrix}
\rho  & & & \\ & 1 & & \\ & & \ddots & \\ & & & 1 
\end{matrix}
\right),
$$
and  
$s^{-1} dx_1 = \rho \thinspace dx_1$, $s^{-1} dx_2 = dx_2$, $\ldots$, 
$s^{-1} dx_n = dx_n$.

Let $I =  \{ I_1, I_2, ..., I_p \} \in \I^p$.
If $I_1 = 1$, then
$$
  \begin{array}{rcl}
     s^{-1} ( \mu_I \hspace{.5em} dx_{I_1} \wedge \cdots \wedge dx_{I_p}) & = &
      s^{-1} \mu_I  \hspace{.5em} s^{-1}
      dx_1 \wedge \cdots \wedge s^{-1}dx_{I_p}\\
      & = & \mu_I \circ s \hspace{.5em} \rho \hspace{.5em} 
      dx_1 \wedge \cdots \wedge dx_{I_p}.
  \end{array}
$$   
If $I_1 \neq 1$, then 
$$
  \begin{array}{rcl}
     s^{-1} ( \mu_I \hspace{.5em} dx_{I_1} \wedge \cdots \wedge dx_{I_p}) & = &
      s^{-1} \mu_I  \hspace{.5em} s^{-1}
      dx_{I_1} \wedge \cdots \wedge s^{-1}dx_{I_p}\\
      & = & \mu_I \circ s \hspace{.5em} 
      dx_{I_1} \wedge \cdots \wedge dx_{I_p}.
  \end{array}
$$
But $\mu$ is \chivar, so
$\rho^{a} \mu = \det(s)^a \mu =\chi^{-1}(s)\mu = s^{-1} \mu$. 
Hence if $I_1 = 1$, then 
  $\rho^a \mu_I = \rho^{\phantom{a}} \mu_I \circ s$,
i.e.
$\rho^{a-1} \mu_I = \mu_I \circ s$.
Thus $x_1^{a-1}$ divides $\mu_I$.  
Similarly, if $I_1 \neq 1$, then 
$\rho^a \mu_I = \mu_I \circ s $
and $x_1^{a}$ divides $\mu_I$.  

\end{proof}

\begin{lemma}
\label{lemma:Qdivides}
  $Q_{\chi}$ divides the exterior product of any two \chivar{}  
  differential forms.
\end{lemma}

\begin{proof}
Let $\mu$ be a \chivar{} $p$-form and $\omega$ be a \chivar{} 
$q$-form.   
 Fix $H \in \A$.  Let $s = s_H$ and $a = a_H(\chi)$.  Assume that
$a \neq 0$.  We show that
$\alpha_H^a$ divides $\mu \wedge \omega$ by choosing coordinates from 
Lemma~\ref{lemma:degx1} in which
$\alpha_H = x_1$. 
Let
$$
\begin{array}{rccl}
    \mu & = & \displaystyle \sum_{I \in \I^p} &\mu_I \thinspace dx_{I_1} 
       \wedge \hdots \wedge dx_{I_p},\\
    \omega & = & \displaystyle \sum_{J \in \I^q} &\omega_J \thinspace dx_{J_1} 
        \wedge \hdots \wedge dx_{J_q}, 
         \text{ and }\\
    \mu \wedge \omega &= & 
         \displaystyle\sum_{K \in \I^{p+q}} &\gamma_K \thinspace dx_{K_1}
        \wedge \hdots \wedge dx_{K_{p+q}}
\end{array}
$$
in these coordinates.  
Then $x_1^a$ divides $\mu_I$ whenever $I_q \neq 1$ and 
$x_1^a$ divides $\omega_J$ whenever $J_1 \neq 1$.  

Hence, for $I \in \I^p$ and $J \in \I^q$, the polynomial $\mu_I \omega_J$
is divisible by $x_1^a$ given that not both $I_1$ and $J_1$ are $1$.  Since
each $\gamma_K$ is either zero or a sum of terms of the form $\pm \mu_I \omega_J$
where the multiindices $I$ and $J$ are disjoint, $x_1^a$ divides
each $\gamma_K$ and hence $\mu \wedge \omega$.  Thus, $\mu \wedge \omega$
is divisible by $\alpha_H^a= \alpha_H^{a_H(\chi)}$.  Since $H$ was arbitrary,
$Q_\chi$ divides $\mu \wedge \omega$. 

\end{proof}

Lemma~\ref{lemma:Qdivides} prompts us to define the following multiplication
in $\Omega^{\chi}$:
For differential forms $\mu$ and $\omega$, define the {\em $\chi$-wedge} of
$\mu$ and $\omega$ as
$$
  \mu \cw \omega := \frac{\mu \wedge \omega}{Q_{\chi}}.
$$
If $\mu$ and $\omega$ are \chivar{} forms, 
$\mu \cw \omega$ is again \chivar.  Thus, Lemma~\ref{lemma:Qdivides}
implies 
\begin{cor}
 \label{cor:closed}
  The $R$-module $\Omega^{\chi}$ is closed under $\chi$-wedging.
\end{cor}

The following proposition gives a condition (similar to
Saito's Criterion) for $n$ $1$-forms to generate $\Omega^{\chi}$.   The
proof is similar to Solomon's original argument \cite{Sol1} that $df_1, \ldots, df_n$ generate 
the module of invariant differential forms.  

\begin{prop}
 \label{prop:generalSaito}
  Let  $\omega_1, \ldots, \omega_n$ be \chivar{} $1$-forms.  The forms
  $\omega_{I_1} \cw \ldots \cw \omega_{I_p}$, for $I \in \I^p$ and  
     $p\geq 0$,
  generate $\Omega^{\chi}$ over $R$ if
  and only if 
  $$
    \omega_1 \cw \ldots \cw \omega_n \doteq Q_{\chi \cdot \det} 
    \thinspace vol.
  $$ 
\end{prop}

\begin{proof}
Assume that $\omega_1 \cw \ldots \cw \omega_n \doteq Q_{\chi \cdot \det} 
    \thinspace vol$.
The $p$-forms
$
    \omega_{I_1} \cw \ldots \cw \omega_{I_p}, I \in \I^p,
$ 
are \chivar{} by Corollary~\ref{cor:closed}. 

Since 
$\omega_1 \cw \ldots \cw \omega_n \neq 0$, 
$\omega_1 \wedge \ldots \wedge \omega_n \neq 0$, and the forms
$\omega_{I_1} \wedge \ldots \wedge \omega_{I_p}$, $I \in \I^p$, are
linearly independent over $F:=\C(x_1, \ldots, x_n)$.  If not, there
exist rational functions $r_I$ with
$$
   0= \displaystyle\sum_{I \in \I^p} r_I \thinspace \omega_{I_1} 
   \wedge \cdots \wedge \omega_{I_p}.
$$
Fix $J \in \I^p$ and $J^c \in \I^{n-p}$. Then
\begin{eqnarray*}
     0 & =&   \left(\displaystyle\sum_{I \in \I^p} r_I \thinspace \omega_{I_1} 
             \wedge \cdots \wedge \omega_{I_p}\right) \wedge 
             \omega_{J_1^c} \wedge \cdots \wedge \omega_{J_{n-p}^c}        \\
       & =&  \vphantom{\left(\displaystyle\sum \right)}
              \pm r_J \thinspace \omega_1 \wedge \cdots \wedge \omega_n,     
 \end{eqnarray*}
and $r_J$ must be zero.  Hence the forms
$$
  \omega_{I_1} \cw \ldots \cw \omega_{I_p}=
   (Q_{\chi})^{1-p} \thinspace \omega_{I_1} \wedge \ldots \wedge \omega_{I_p},
    \quad I \in \I^p,
$$ 
are also linearly independent over $F$, and thus span 
$$ 
    \Omega^{p}(V)  := \bigoplus_{I \in \I^p}
    F dx_{I_1} \wedge \ldots \wedge dx_{I_p}     
$$
since $\Omega^{p}(V)$ has dimension $\binom{n}{p}$.

Choose an arbitrary \chivar{} $p$-form $\mu$.  Then there exist rational 
functions $t_I$ $\in F$ with 
$$
    \mu = \displaystyle \sum_{I \in \I^p} t_I\thinspace
            \omega_{I_1} \cw \cdots \cw \omega_{I_p}.
$$
Fix $J \in \I^p$ and its complementary index $J^c$. We will
show that $t_J \in R$.

By Corollary~\ref{cor:closed}, 
the $n$-form
$(\omega_{J_1^c} \cw \cdots \cw \omega_{J_{n-p}^c}) \cw \mu$
is \chivar{}.  Thus by Equation (\ref{invn-forms}) above, 
there exists a polynomial $f \in R$ with
$$
  \left( \omega_{J_1^c} \cw \cdots \cw \omega_{J_{n-p}^c}\right) \cw \mu
   =f \thinspace  Q_{\chi \cdot \det} \thinspace vol.
$$
On the other hand,
\begin{align*}
    \left( \omega_{J_1^c} \cw \right. 
     \cdots \cw & \left. \omega_{J_{n-p}^c}\right) \cw\mu \\  
    & =    
    \left( \omega_{J_1^c} \cw \cdots \cw \omega_{J_{n-p}^c} \right) \cw
    \displaystyle\sum_{I \in \I^p} t_I \thinspace 
                \omega_{I_1} \cw \cdots \cw \omega_{I_p}\\
    & = 
    \left( Q_{\chi}^{1-n} \right) 
    \left( \omega_{J_1^c} \wedge  
             \cdots \wedge \omega_{J_{n-p}^c}\right) \wedge
           \displaystyle \sum_{I \in \I^p} t_I \thinspace 
                \omega_{I_1} \wedge \cdots \wedge \omega_{I_p}\\
    & = 
    \left( Q_{\chi}^{1-n} \right) 
         \vphantom{\displaystyle\sum_{i \in \I^p}}\pm\displaystyle
         t_J \thinspace 
          \omega_{1} \wedge \cdots \wedge \omega_{n}\\
    & =  \vphantom{\displaystyle\sum_{\I^p}}
       \pm t_J \thinspace 
         \omega_{1} \cw \cdots \cw \omega_{n} \\  
    & \doteq
       \vphantom{\displaystyle\sum_{\I^p}}
       \pm t_J 
        \thinspace Q_{\chi \cdot \det} \thinspace vol.\\
\end{align*}

\vspace{-6ex} 
\noindent
Thus
$f \thinspace Q_{\chi \cdot \det}\thinspace 
  \doteq \pm t_J 
  \thinspace Q_{\chi \cdot \det}$.
Hence, $t_J \in R$.  Since $J$ was arbitrary, $\mu$ is
in the $R$-span of 
$\{\omega_{I_1} \cw \ldots \cw \omega_{I_p}, \hspace{.5ex} I_p \in \I^p \}.$

The converse follows from Equation (\ref{invn-forms}) above.

\end{proof}
\section{Condition satisfied}

Since $\Omega^p$ has rank $\binom{n}{p}$, the 
$R$-module $(\Omega^p)^{\chi}$ is also free of rank $\binom{n}{p}$ (this
follows from Lemma 6.45 of \cite{OrTer}, p.\hspace{1ex}232).
We will show that the generators of $(\Omega^1)^{\chi}$
satisfy the condition given in Proposition~\ref{prop:generalSaito}, 
but first we must 
gather some preliminary facts.  

We recall some results about invariant vector fields.  There exist $n$
invariant vector fields, called {\em basic derivations}, that generate the
module of invariant vector fields over $R$ (see \cite{OrTer}, Section 6.3).  Using
Saito's Criterion, H. Terao showed that the coefficient matrix of the basic
derivations has determinate $Q_{\det^{-1}}$ upto a nonzero scaler multiple
(see \cite{OrTer}, 
p.\hspace{1ex}238). 
Using the minors of this coefficient matrix, we 
construct $(\det^{-1})\hspace{.2ex}$-invariant  $1$-forms, 
$\mu_1, \ldots, \mu_n$,
that satisfy
\begin{eqnarray*}
      \mu_1\wedge\cdots\wedge \mu_n                                   
   &=& Q_{\det^{-1}}^{n-1} \thinspace vol. 
\end{eqnarray*}
The forms $\mu_1, \ldots,
\mu_n$ thus generate 
$\Omega^{\det^{-1}}$ over $R$ by Proposition~\ref{prop:generalSaito}.  We will use these forms to give an
argument for arbitrary $\chi$.  

We also note the relationship between $Q_{\chi\cdot\det}$ and 
$Q_{\chi}$:
Fix $H \in \A$ with $a_H(\chi) \neq 0$.  The exponent
$a_H(\chi\cdot\det)$ is the least nonnegative integer 
satisfying 
\begin{align*}
  \det(s_H)^{-a_H(\chi\cdot\det)} 
              & = ( \chi\cdot\det)(s_H)                       \\
              & = \chi(s_H)\thinspace \det(s_H)              \\  
              & = \det(s_H)^{-a_H(\chi)} \thinspace \det(s_H)        \\
              & = \det(s_H)^{-( a_H(\chi)-1)}.   
\end{align*}
Hence, $a_H(\chi\cdot\det)  = a_H(\chi)-1$.   
Now fix $H \in \A$ with $a_H(\chi) = 0$. Then
\begin{align*}
  \det(s_H)^{-a_H(\chi\cdot\det)} 
              & = ( \chi\cdot\det)(s_H)                      \\
              & = \chi(s_H)\thinspace \det(s_H)              \\  
              & = \det(s_H)        \\
              & = \det(s_H)^{-(o(s_H)-1)},
\end{align*}
and $a_H(\chi\cdot\det)  = o(s_H)-1$.   
Thus,
\begin{eqnarray*}
   Q_{\chi\cdot\det} &=& \prod_{H\in\A}\alpha_H^{{a_H}(\chi\cdot\det)}\\
                     &=& \underset{\chi(s_H) \neq 1}
                                 {\prod_{H\in\A}}         
                          \alpha_H^{a_H(\chi)-1} \thinspace
                          \underset{\chi(s_H) = 1}
                                 {\prod_{H\in\A}}         
                          \alpha_H^{o(s_H)-1}.                          
\end{eqnarray*}


\begin{prop}
  \label{prop:criterion-satisfied}
  If $\omega_1, \ldots, \omega_n$ generate $(\Omega^1)^{\chi}$ 
  over $R$, then
  $$
      \omega_1 \cw \ldots \cw \omega_n \doteq Q_{\chi \cdot \det} 
      \thinspace vol.
  $$ 
\end{prop}

\begin{proof}
Let $M$ be the coefficient matrix of $\omega_1, \ldots, \omega_n$,
i.e. $\omega_1 \wedge \ldots \wedge \omega_n = \det M \thinspace vol$. 
Suppose that $\det M=0$.  Then one row of $M$ is a linear
combination of the other rows over $F=\C(x_1, \ldots, x_n)$.  Multiplying
by a least common multiple yields a relation over $S$:
$\sum_{i=1}^n s_i \thinspace \omega_i = 0$.
To get a relation over $R$, apply a group element $g$, multiply by $\chi^{-1}(g)$,
and then sum over $G$:
\begin{eqnarray*}
   0&=& \displaystyle\sum_{g \in G} \displaystyle\sum_{i=1}^n
         \chi^{-1}(g) \thinspace g s_i \thinspace g\omega_i \\
    &=& 
   \displaystyle\sum_{i=1}^n \displaystyle\sum_{g \in G} 
         \chi^{-1}(g) \thinspace g s_i \thinspace \chi(g) \thinspace \omega_i \\ 
    &=&
    \displaystyle\sum_{i=1}^n \left(\displaystyle\sum_{g \in G} 
                  g s_i \right) \omega_i .
\end{eqnarray*}
This contradicts the fact that $(\Omega^1)^{\chi}$ is free over $R$ 
with basis $\omega_1, \ldots, \omega_n$.
Hence, $\det M \neq 0$.

By Corollary~\ref{cor:closed}, 
$\omega_1 \cw \cdots\cw \omega_n$ is a \chivar{} $n$-form.  Thus (from 
Equation (\ref{invn-forms}))
there exists a nonzero $f \in R$ with
$$
   (Q_{\chi})^{1-n} \thinspace \det M \thinspace vol 
    =  (Q_{\chi})^{1-n} \thinspace \omega_1 \wedge\ldots\wedge \omega_n
    =  \omega_1 \cw \cdots\cw \omega_n 
    =  f \thinspace Q_{\chi\cdot\det} \thinspace vol.
$$
Hence, $\det M=f \thinspace Q_{\chi\cdot\det} \thinspace (Q_{\chi})^{n-1}$.

We show that $f$ is constant by finding two polynomials
that share no factors, yet are each divisible by $f$. 
Since each 
$df_i$ is invariant, each $Q_{\chi} \thinspace df_i$ is \chivar{} and 
hence a combination
of $\omega_1, \ldots, \omega_n$ over $R$. There exists a matrix of 
coefficients, 
$N$, with entries in $S$, such that
\begin{align*}
   Q_{\chi} df_1 \wedge \cdots\wedge Q_{\chi}df_n 
   &= \det M \thinspace \det N \thinspace  vol \\
   &=  f \thinspace Q_{\chi\cdot\det} \thinspace (Q_{\chi})^{n-1}
      \thinspace \det N \thinspace vol.
\end{align*}
But, $df_1 \wedge\cdots\wedge  df_n \doteq Q_{\det}$, so
$$
    Q_{\chi} df_1 \wedge\cdots\wedge Q_{\chi} df_n 
   \doteq (Q_{\chi})^n \thinspace Q_{\det} \thinspace vol.
$$
Hence, 
$$
  f \thinspace Q_{\chi\cdot\det} \thinspace \det(N) 
  \doteq Q_{\chi} \thinspace Q_{\det}
$$
and since 
$\det N \in S$, $f$ divides 
$Q_{\chi} \thinspace Q_{\det} \thinspace (Q_{\chi\cdot\det})^{-1}$.

Since each $\mu_i$ (introduced above) is 
$(\det^{-1})\hspace{.2ex}$-invariant, each 
$Q_{\chi\cdot\det} \thinspace \mu_i$ is \chivar{}, and thus a $R$-combination
of $\omega_1, \ldots, \omega_n$.  There exists a matrix
of coefficients, $N'$, with coefficients in S, such that 
\begin{align*}
  Q_{\chi\cdot\det} \mu_1 \wedge \cdots\wedge Q_{\chi\cdot\det} \mu_n 
   &= \det M  \thinspace \det N'\thinspace  vol \\
   &=  f \thinspace Q_{\chi\cdot\det} \thinspace (Q_{\chi})^{n-1}
      \thinspace \det N' \thinspace vol.
\end{align*}
But we choose the $\mu_i$ so that
$$
   Q_{\chi\cdot\det} \mu_1 \wedge\cdots\wedge Q_{\chi\cdot\det} \mu_n 
   = (Q_{\chi\cdot\det})^n\thinspace (Q_{\det^{-1}})^{n-1}\thinspace vol.
$$
Hence, 
$$
  f \thinspace Q_{\chi\cdot\det} \thinspace 
   \left(Q_{\chi}\right)^{n-1} \thinspace \det N'
  = \left(Q_{\chi\cdot\det}\right)^n (Q_{\det^{-1}})^{n-1}
$$
and since 
$\det N' \in S$, 
$\left(Q_{\chi\cdot\det} Q_{\det^{-1}}\right)^{n-1}
     \thinspace\left(Q_{\chi}\right)^{1-n}$
is divisible by $f$.

We show that the two polynomials
$$Q_{\chi} \thinspace Q_{\det} \thinspace (Q_{\chi\cdot\det})^{-1}
   \text{ and }
\left(Q_{\chi\cdot\det} Q_{\det^{-1}}\right)^{n-1}
     \thinspace\left(Q_{\chi}\right)^{1-n}
$$
have no common factors
by writing them both in terms of the $\alpha_H$.  We 
expand the factors:
\begin{eqnarray*}
  Q_{\chi} &=& \underset{\chi(s_H) \neq 1}
                          {\prod_{H\in\A}}         
                          \alpha_H^{a_H(\chi)} \thinspace ,  \\           
  Q_{\det} &=&  \underset{\chi(s_H) \neq 1}
                          {\prod_{H\in\A}}         
                          \alpha_H^{o(s_H)-1} 
                          \thinspace
                          \underset{\chi(s_H) = 1}
                                 {\prod_{H\in\A}}         
                          \alpha_H^{o(s_H)-1} \thinspace ,    \\          
  Q_{\chi\cdot\det}  &=&
                          \underset{\chi(s_H) \neq 1}
                          {\prod_{H\in\A}}         
                          \alpha_H^{a_H(\chi)-1} \thinspace
                          \underset{\chi(s_H) = 1}
                                 {\prod_{H\in\A}}         
                          \alpha_H^{o(s_H)-1}     \thinspace    ,      \\
  Q_{\det^{-1}}  &=&
                          \underset{\chi(s_H) \neq 1}
                          {\prod_{H\in\A}}         
                          \alpha_H^{\phantom{a_H(\chi)-1}} 
                          \thinspace
                          \underset{\chi(s_H) = 1}
                                 {\prod_{H\in\A}}         
                          \alpha_H  \thinspace.
\end{eqnarray*}
The first polynomial, 
$Q_{\chi} \thinspace Q_{\det} \thinspace (Q_{\chi\cdot\det})^{-1}$, 
simplifies to
$$
\underset{\chi(s_H) \neq 1}
                          {\prod_{H\in\A}}         
                          \alpha_H^{o(s_H)} \thinspace,  
$$
and the second polynomial, 
$\left(Q_{\chi\cdot\det}\thinspace Q_{\det^{-1}}\right)^{n-1}
\thinspace\left(Q_{\chi}\right)^{1-n}$, simplifies to
$$
   \left( \underset{\chi(s_H) = 1}
                                 {\prod_{H\in\A}}         
                          \alpha_H^{o(s_H)}                         
    \right)^{n-1} .
$$

Since $f$ divides both polynomials, $f$ must be constant.  
Thus, 
$\omega_1, \ldots, \omega_n$ satisfy the criterion of
Proposition~\ref{prop:generalSaito}.

\end{proof}


\begin{cor}
  \label{cor:free-like}
  There exist $n$ $1$-forms $\omega_1, \ldots, \omega_n$ such that  
  $\Omega^{\chi}$ is generated over $R$ by the forms
  $\omega_{I_1} \cw \ldots \cw \omega_{I_p}, \hspace{.5ex}I \in \I^p, 
  \hspace{.5ex} p\geq 0$.  Thus  $\Omega^{\chi}$ has the
  structure of an exterior algebra.  
  
\end{cor}

\section{Example:  $G_{26}$}
For an example, let us take a three dimensional complex reflection group,
$G_{26}$.  This group is the symmetry group of a regular complex
polyhedron,
and is number 26 in Shephard and Todd's enumeration of
finite irreducible unitary groups generated by reflections
\cite{Todd-She}.  The group $G_{26}$ consists of 1,296 complex $3 \times
3$ matrices and is generated by reflections of order two and three.  
The associated
collineation group (which results from moding out by the 
scaler matrices) is the
Hessian group of order 216.  

The group is 
generated by the matrices
\begin{small}
$$
\left(
\begin{array}{ccc}
1  &    0   &  0 \\
0  &    0   &  1\\
0  &    1   &  0
\end{array}
\right),
\hspace{2ex}
\left(
\begin{array}{ccc}
1  &    0   &  0 \\
0  &    1   &  0\\
0  &    0   &  \alpha^2
\end{array}
\right),
\hspace{2ex} \text{\normalsize and} \hspace{2ex}
\frac{i}{\sqrt{3}}
 \left(
\begin{array}{ccc}
\alpha  &  \alpha^2 & \alpha^2  \\
\alpha^2  &  \alpha & \alpha^2  \\
\alpha^2  &  \alpha^2 & \alpha  
\end{array}
\right),
$$
\end{small}
where $\alpha$ is a primitive cube root of unity.   

The character table for this group reveals six multiplicative characters,
each a power of the determinate character.  
Choose $\chi = \det^3$. Note that
$$
Q_{\det^3} = (x^3-y^3)(x^3-z^3)(y^3-z^3), 
$$
and
$$
Q_{\det^4} =  
x^2 y^2 z^2 (x^9 + 3 x^6 (y^3 + z^3) + (y^3 + z^3)^3 + 
     3 x^3 (y^6 - 7 y^3 z^3 + xz^6))^2.
$$

The following $1$-forms are $\det^3$-invariant:   
\begin{eqnarray*}
 \omega_1 & = &
 \hphantom{-}x^2(y - z)(y^2 + y z + z^2)
 (\hphantom{-}2x^3 - y^3 - z^3) \thinspace dx \\
  &&-y^2(x - z)(x^2 + x z + z^2)(-x^3 + 2 y^3 - z^3)\thinspace dy \\ 
  &&  - z^2(x - y)(x^2 + x y + y^2)(\hphantom{-}x^3 + y^3 - 2z^3)
   \thinspace dz, \\
\end{eqnarray*}

\vspace{-7ex}

\begin{eqnarray*}
  \omega_2 & = & 
   x^2 (x^3 - y^3)(x^3 - z^3)(y^3 - z^3)(\hphantom{-}x^3 - 5y^3 -5 z^3)
   \thinspace dx\\
   && y^2 (x^3 - y^3)(x^3 - z^3)(y^3 - z^3)(-5x^3 + y^3 - 5z^3)\thinspace dy\\
   && z^2 (x^3 - y^3)(x^3 - z^3)(y^3 - z^3)(-5x^3 - 5y^3 + z^3)\thinspace dz,\\
\end{eqnarray*}

\vspace{-7ex}

\begin{eqnarray*}
  \omega_3 & = &
   x^2(x^3 - y^3)(x^3 - z^3)(y^3 - z^3)
   (x^9 + 3y^9 + 61y^6z^3 +  61y^3z^6 
    + 3z^9 \\
  &&  \hspace{7ex} 
      + 9x^6(y^3 + z^3) + x^3(-13y^6 + 122y^3z^3 - 13z^6))\thinspace dx + \\
  &&  y^2(x^3 - y^3)(x^3 - z^3)(y^3 - z^3)
    (3x^9 + y^9 + 9y^6z^3  - 13y^3z^6 + 3z^9 \\
  && \hspace{7ex} 
     + x^6(-13y^3 + 61z^3) + x^3(9y^6 + 122y^3z^3 + 61z^6))\thinspace dy +\\ 
  && z^2 (x^3 - y^3)(x^3 - z^3)(y^3 - z^3)
    (3x^9 + 3y^9 - 13y^6z^3 + 9y^3z^6 + z^9 \\
  && \hspace{7ex}
      + x^6(61y^3 - 13z^3) + x^3(61y^6 + 122y^3z^3 + 9z^6))\thinspace dz.
\end{eqnarray*}
The polynomial $Q_{\det^3}$ divides 
$\omega_1 \wedge \omega_2$, $\omega_2 \wedge \omega_3$,
and $\omega_1 \wedge \omega_3$.  The 
determinate of the coefficient matrix of $\omega_1$, $\omega_2$,
and $\omega_3$ is $(-16) Q_{\det^4} \thinspace Q_{\det^3}^2$,
hence $\omega_1$, $\omega_2$,
and $\omega_3$ $\chi$-wedge to a multiple of $Q_{\chi\cdot\det} = Q_{\det^4}$.  
Proposition~\ref{prop:generalSaito} then implies that
$\omega_1$, $\omega_2$, and $\omega_3$ generate the entire module
of $\det^3$-invariants over the ring of invariants via $\det^3$-wedging.

\section{Logarithmic forms}
We have so far only discussed regular differential forms; we now consider
rational differential forms. 
The $S$-module of {\em logarithmic $p$-forms with poles along $\A$} 
(see also \cite{OrTer}, p. 124) is 
defined as
$$
  \Omega^p(\A) := \{
    \frac{\omega}{Q_{\det^{-1}}} : \omega \in \Omega^p \text{ and }
     \omega \wedge d \alpha_H \in \alpha_H  \thinspace \Omega^{p+1} 
     \text{ for all } H \in \A \}.
$$
Ziegler \cite{Ziegler} extends this definition to 
{\em multiarrangements of hyperplanes}, hyperplane arrangements in which each
hyperplane has a positive integer multiplicity.  We apply  
his definitions to our context of reflection groups and semiinvariants: 
Let $\A_{\chi}$ be the multiarrangement consisting of 
hyperplanes $H \in \A$ each with multiplicity
$\alpha_H(\chi)$, i.e. the multiarrangement defined by $Q_\chi$.
We define (as in \cite{Ziegler}) the module of 
{\em logarithmic $p$-forms} of $\A_\chi$:
$$
  \Omega^p(\A_\chi) := \{
    \frac{\omega}{Q_\chi} : \omega \in \Omega^p \text{ and }
     \omega \wedge d \alpha_H \in \alpha_H^{a_H(\chi)} \thinspace \Omega^{p+1} 
     \text{ for all } H \in \A \}.
$$
Let
$$
\Omega(\A_\chi)  :=  \bigoplus_{p \geq 0} \Omega^p(\A_\chi).
$$

\begin{cor}
  $$
     \Omega^\chi \subset Q_\chi \thinspace \Omega(\A_\chi).
  $$
\end{cor}
\begin{proof}
Choose $\omega$ in $(\Omega^p)^\chi$ and fix $H \in \A$.  
Using Lemma~\ref{lemma:degx1}, 
choose coordinates in which $x_1 = \alpha_H$, 
$\omega = \sum_{I \in \I^p} \omega_I \thinspace dx_{I_1} \wedge\hdots
\wedge dx_{I_p}$,  and
$x_1^{a_H(\chi)}$ divides $\omega_I$ if $1 \notin I$.  
Then $d\alpha_H = dx_1$, and 
$\omega \wedge d \alpha_H = \omega \wedge dx_1= \sum_{I, 1 \notin I} \omega_I
\wedge dx_1$, 
which is divisible by $x_1^{a_H(\chi)}$.  Hence, 
$ \omega \wedge d \alpha_H \in \alpha_H^{a_H(\chi)} \thinspace 
\Omega^{p+1}$. As $H$ was arbitrary, $\frac{\omega}{Q_\chi} \in \Omega(\A_\chi)$.

\end{proof}

This relationship is stronger when $\chi = \det^{-1}$.  In this case, 
the forms that generate $\Omega^\chi$ via $\chi$-wedging over $R$ also
generate $\Omega(\A)$ over $S$ (see \cite{Terao-Shepler} for an 
independent proof).

On a similar note, we have
\begin{prop}
   $\Omega(\A_\chi)$ is closed under the exterior product.
\end{prop}
\begin{proof}
Let $\omega / Q_\chi$ and $\mu / Q_\chi$ be in $\Omega(\A_{\chi})$. 
Fix $H$ in $\A$ and let $a_H(\chi)= a$.
Choose coordinates such that $x_1 =
\alpha_H$, and write
$\omega = \sum_{I \in \I^p} \omega_I \thinspace dx_{I_1} \wedge\hdots
\wedge dx_{I_p}$ and  $\mu = \sum_{J \in \I^q} \mu_J \thinspace dx_{J_1} \wedge\hdots
\wedge dx_{J_q}$ in these coordinates.
Since $\omega \wedge dx_1 = \omega \wedge d\alpha_H \in
\alpha_H^a \Omega = x_1^a \Omega$, 
$\omega_I$ is divisible by $x_1^a$ as long as $1 \not\in I$.  Similarly,
$\mu_J$ is divisible by $x_1^a$ whenever $1 \not\in J$.  
As in the proof of Lemma~\ref{lemma:Qdivides}, it follows that $Q_\chi$ divides 
$\omega \wedge \mu$. 
Whenever
$1\not\in I$ {\em and} $1\not\in J$, $x_1^{2a}$ divides $\omega_I \mu_J$, and
thus 
$$
  \frac{\omega \wedge \mu}{Q_\chi} \wedge dx_1
$$ 
is also divisible by $x_1^a$.  Hence $\alpha_H^a$ divides
$(1/Q_\chi)\thinspace \omega
\wedge\mu\wedge d \alpha_H$, and as $H$ was arbitrary, 
$(\omega /Q_\chi) \wedge (\mu /Q_\chi)$ is in
$\Omega(\A_\chi)$.

\end{proof}

\section{Remarks}

Analogous results hold for vector fields, or {\em derivations}.
Let $\Upsilon^\chi$ be the module of \chivar s
in the exterior algebra of derivations. 
Because the group action differs here, 
Lemma~\ref{lemma:degx1} is slightly different,
with $a+1$ taking the place of $a-1$ when
$I_1 =1$.  
The case where $I_1 \neq 1$ is the same as in the original lemma, and hence 
$Q_\chi$ also divides the exterior product of two
elements in $\Upsilon^\chi$ (the proof is analogous to the
case of $\Omega^{\chi}$).  The criterion for $n$ derivations
to generate $\Upsilon^\chi$ via $\chi$-wedging is also slightly different:
they must $\chi$-wedge to 
$Q_{\chi\cdot\det^{-1}} \thinspace \frac{\del}{\del x_1}
\wedge\ldots\wedge  \frac{\del}{\del x_n}$  
instead of $Q_{\chi\cdot\det} \thinspace dx_1 \wedge\ldots\wedge dx_n$.  This follows from the fact
that $dx_1 \wedge\ldots\wedge dx_n$ is $(\det^{-1})$-invariant while
$\frac{\del}{\del x_1} \wedge\ldots\wedge  \frac{\del}{\del x_n}$
is $\det$-invariant.  Finally, we note that 
the correspondence between differential $p$-forms
(in $\Omega^p$) and $(n-p)$-forms in $\Upsilon$ (the exterior algebra 
of derivations) induces a module isomorphism between
$\Omega^\chi$ and $\Upsilon^{\chi\cdot\det}$.

\section{Acknowledgments}
The author is grateful to Peter Doyle, Hiroaki Terao, and Nolan Wallach
for their helpful comments.



\end{document}